\documentclass{amsart}
\usepackage{amssymb}
\usepackage{amsmath}
\newtheorem{theorem}{Theorem}[section]

\newtheorem{proposition}[theorem]{Proposition}

\newcommand{\mpn}{\medskip\par\noindent}
\newcommand{\pn}{\par\noindent}

\newcommand{\bpn}{\bigskip\par\noindent}
\theoremstyle{definition}

\numberwithin{equation}{section}

\begin{document}
\newcommand{\Mod}[1]{\,(\text{\mbox{\rm mod}}\;#1)}
\title[  Identities for Bernoulli polynomials involving Chebyshev polynomials  ]{Some identities for Bernoulli polynomials involving Chebyshev polynomials}
\author{DAE SAN KIM, TAEKYUN KIM AND Sang-Hun Lee}
\begin{abstract}
In this paper  we derive  some new and interesting identities for     Bernoulli, Euler  and  Hermite polynomials associated with Chebyshev polynomials.
\end{abstract}
\maketitle
\vspace{3mm}

\def\ord{\text{ord}_p}
\par\bigskip\noindent
\section{Introduction }
The Bernoulli number are defined by the generating function to be
\begin{equation}\tag{1}
\begin{split}
\frac{t}{e^t-1}=e^{Bt}=\sum_{n=0}^{\infty}\frac{B_n}{n!}t^n, \quad (\text{see [3,13,14]}),
\end{split}
\end{equation}
with the usual convention about replacing $B^n$ by $B_n$.

As is well known, the Bernoulli polynomials are given by
\begin{equation}\tag{2}
\begin{split}
B_n(x)=(B+x)^n=\sum_{l=0}^{n}\binom{n}{l}B_{n-l}x^l, \quad (\text{see [1-8]}).
\end{split}
\end{equation}
From (1), we note that the recurrence relation for the Bernoulli numbers is given by
\begin{equation*}
\begin{split}
B_0=1,   \quad    (B+1)^n-B_n=\delta_{1,n}, \quad (\text{see [6-8]}),
\end{split}
\end{equation*}
where $\delta_{m,n}$ is the Kronecker symbol.

By (2), we get
\begin{equation}\tag{3}
\begin{split}
\frac{dB_n(x)}{dx}=  n \sum_{l=0}^{n-1}\binom{n-1}{l}B_{n-1-l}x^l=nB_{n-1 }(x).
\end{split}
\end{equation}
Thus, by (3), we see that
\begin{equation}\tag{4}
\begin{split}
\int B_n(x) dx = \frac{B_{n+1}(x)}{n+1}+C, \quad (\text{see [3]}),
\end{split}
\end{equation}
where $C$ is a some constant.

The Euler polynomials are defined by the generating function to be
\begin{equation}\tag{5}
\begin{split}
\frac{2}{e^t+1}e^{xt}=e^{E(x)t}=\sum_{n=0}^{\infty}E_n(x)\frac{t^n}{n!},
\end{split}
\end{equation}
with the usual convention about replacing $E^n(x)$ by $E_n(x)$, (see [1,2,4,10,11]).

In the special case, $x=0$, $E_n(0)=E_n$ are   called the $n$-th Euler numbers.

It is well known $[6,15]$ that Hermite polynomials are given by the generating function to be
\begin{equation}\tag{6}
\begin{split}
e^{2xt-t^2}=e^{H(x)t}=\sum_{n=0}^{\infty}H_n(x)\frac{t^n}{n!},
\end{split}
\end{equation}
with the usual convention about replacing $H^n(x)$ by $H_n(x)$.

From (6),we have
\begin{equation}\tag{7}
\begin{split}
\frac{dH_n(x)}{dx}= 2 nH_{n-1}(x), \quad H_n(x)=(-1)^nH_n(-x).
\end{split}
\end{equation}
By (1) and (2), we easily get
\begin{equation}\tag{8}
\begin{split}
B_n(x)=\sum_{\substack{k=0 \\ k\neq 1} }^{n }\binom{n}{k}E_{n-k}(x), \quad (\text{see [1-15]}),
\end{split}
\end{equation}
\begin{equation}\tag{9}
\begin{split}
E_n(x)=-2\sum_{l=0 }^{n }\binom{n}{l}\frac{E_{l+1}}{l+1}E_{n-l}(x),
\end{split}
\end{equation}
and
\begin{equation}\tag{10}
\begin{split}
x^n=\frac{1}{n+1}\bigl(B_{n+1}(x+1)-B_{n+1}(x )\bigl)=\frac{1}{n+1}\sum_{l=0 }^{n }\binom{n+1}{l}B_{l}(x ).
\end{split}
\end{equation}
The Chebyshev polynomial $T_n(x)$ of the first kind is a polynomial in $x$ of degree $n$, defined by the relation
\begin{equation}\tag{11}
\begin{split}
T_n(x)=\cos{n\theta},\quad\text{when } x=\cos{ \theta}, \quad (\text{see [9]}).
\end{split}
\end{equation}
If the range of the variable $x$ is the interval $[-1,1]$, then the range of the corresponding variable $\theta$
can be taken as $[0,\pi]$. It is known that $\cos{n\theta}$ is a polynomial of degree $n$ in $\cos{ \theta}$, and
indeed we are familiar with elementary formulas $\cos{3\theta}=4\cos^3{ \theta}-3\cos{ \theta}$, $\cos{4 \theta}=8\cos^4{ \theta}-8\cos^2{ \theta}+1$,
$\cdots$.

Thus, by (11), we get
\begin{equation*}
\begin{split}
&T_0(x)=1, \quad T_1(x)=x, \quad T_2(x)=2x^2-1, \quad T_3(x)=4x^3-3x,
\\& T_4(x)=8x^4-8x^2+1,\cdots.
\end{split}
\end{equation*}
The Chebyshev polynomial  $U_n(x)$ of the second kind is a polynomial of degree $n$ in $x$ defined by
\begin{equation}\tag{12}
\begin{split}
U_n(x)=\sin{(n+1)\theta}/\sin{\theta},\quad\text{when } x=\cos{ \theta}, \quad (\text{see [9]}).
\end{split}
\end{equation}
Thus, from (12), we have
\begin{equation*}
\begin{split}
U_0(x)=1, \quad U_1(x)=2x, \quad U_2(x)=4x^2-1, \quad U_3(x)=8x^3-4x,\cdots.
 \end{split}
\end{equation*}
By (11), we see that $T_n(x)$ is a polynomial of degree $n$ with integral coefficients and the leading coefficient $2^{n-1}~(n\geq1)$
and $1~(n=0)$. It is not difficult to show that $U_n(x)$ is a polynomial of degree $n$  with integral coefficients and the leading coefficient
$2^n~(n\geq0)$.
$T_n(x)$ is a solution of $(1-x^2)y^{\prime\prime}-xy^{\prime}+n^2y=0$ and $U_n(x)$ is a solution of $(1-x^2)y^{\prime\prime}-3xy^{\prime}+n(n+2)y=0$.
It is well known [9] that the generating functions of $T_n(x)$ and $U_n(x)$ are given by
\begin{equation}\tag{13}
\begin{split}
\frac{1-xt}{1-2xt+t^2}=\sum_{n=0}^{\infty}T_n(x)t^n,
\end{split}
\end{equation}
and
\begin{equation}\tag{14}
\begin{split}
\frac{1 }{1-2xt+t^2}=\sum_{n=0}^{\infty}U_n(x)t^n, \quad\text{for } |x|\leq1, ~~|t|<1.
\end{split}
\end{equation}
From (11) and (12), we have
\begin{equation}
\tag{15}
\begin{split}
\int_{-1}^{ 1}\frac{T_n(x)T_m(x)}{\sqrt{1-x^2}}dx= \left\{
\begin{array}{cl}
0, & \text{ if} ~~n\neq m  \\
\frac{\pi}{2}, & \text{ if} ~~n=m>0 \\
\pi, & \text{ if} ~~n=m=0 \\
\end{array}\right.,
\end{split}
\end{equation}
and
\begin{equation}
\tag{16}
\begin{split}
\int_{-1}^{ 1}(1-x^2)^{1/2}U_n(x)U_m(x)dx=\frac{\pi}{2}\delta_{n,m}, \quad (\text{see [9]}).
\end{split}
\end{equation}
The equations (15) and (16) are used to derive our main result in this paper.

The Rodrigues' formulae for $T_n(x)$ and $U_n(x)$ are known as follows:
\begin{equation}
\tag{17}
\begin{split}
 T_n(x)=\frac{(-1)^n2^nn!}{(2n)!}(1-x^2)^{1/2} \biggl(\frac{d^n}{dx^n}(1-x^2)^{n-1/2}\biggl),
\end{split}
\end{equation}
and
\begin{equation}
\tag{18}
\begin{split}
 U_n(x)=\frac{(-1)^n2^n(n+1)!}{(2n+1)!}(1-x^2)^{-1/2} \biggl(\frac{d^n}{dx^n}(1-x^2)^{n+1/2}\biggl).
\end{split}
\end{equation}
The equations (17) and (18) are also used to derive our result related to orthogonality of Chebyshev polynomials.

From (11) and (12), we can easily derive the following equations (19) and (20):
\begin{equation}
\tag{19}
\begin{split}
 T_n(x)=\frac{(x+\sqrt{x^2-1})^n+(x-\sqrt{x^2-1})^n}{2},
\end{split}
\end{equation}
and
\begin{equation}
\tag{20}
\begin{split}
 U_n(x)=\frac{(x+\sqrt{x^2-1})^{n+1}-(x-\sqrt{x^2-1})^{n+1}}{2\sqrt{x^2-1}}.
\end{split}
\end{equation}
By the definitions of $T_n(x)$ and $U_n(x)$, we easily get
\begin{equation}
\tag{21}
\begin{split}
\frac{dT_n(x)}{dx}=nU_{n-1}(x), \quad \frac{dU_n(x)}{dx}=\frac{(n+1)T_{n+1}(x)-xU_n(x)}{x^2-1}.
\end{split}
\end{equation}
From (21), we have
\begin{equation}
\tag{22}
\begin{split}
\int U_n(x)dx=\frac{T_{n+1}(x)}{n+1}, \quad \int T_n(x)dx=\frac{nT_{n+1}(x)}{n^2-1}-\frac{xT_{n }(x)}{n -1}.
\end{split}
\end{equation}
In this paper we derive some new and interesting identities for Bernoulli, Euler and Hermite polynomials arising from the orthogonality
of the Chebyshev polynomials for the inner product space with weighted inner product.
\par\bigskip\noindent
\section{ Some identities for Bernoulli, Euler and Hermite polynomials involving Chebyshev polynomials }

Let $\mathbf{P}_{n}=\{p(x)\in  \mathbb{Q}[x] \mid \text{deg}~p(x)\leq n \}$. Then $\mathbf{P}_{n}$ is an inner product space with
the weighted  inner product
 $$\langle p(x),~q(x) \rangle =\int_{-1}^{1}\frac{p(x)q(x)}{\sqrt{1-x^2}}dx,\quad  \text{where } ~~   p(x), q (x) \in\mathbf{P}_{n}.$$
 From (15), we note that  $\{T_{0} (x), T_{1} (x), \cdots ,T_{n} (x) \}$ is
 an orthogonal basis for   $\mathbf{P}_{n}$.
 \\
 Let us assume  $p(x) \in\mathbf{P}_{n}$. Then $p(x)$ is generated by $\{T_{0} (x), T_{1} (x), \cdots ,T_{n} (x) \}$ to be
\begin{equation}
\tag{23}
\begin{split}
p(x)=\sum_{k=0}^{n}C_kT_k(x).
\end{split}
\end{equation}
By (15) and (23), we get
\begin{equation}
\tag{24}
\begin{split}
&C_k=\frac{\delta_{k}}{\pi}\int_{-1}^{1}\frac{T_k(x)p(x) }{\sqrt{1-x^2}}dx=\frac{\delta_{k}}{\pi}\frac{(-1)^k2^kk!}{(2k)!}\int_{-1}^{1}
\biggl(\frac{d^k}{d  x^k }(1-x^2)^{k-1/2}\biggl)p(x)dx,\\
&   \text{where } \delta_{k}= \left\{
\begin{array}{cl}
1, & \text{ if} ~~k=0  \\
2, & \text{ if} ~~k>0. \\
 \end{array}\right.
\end{split}
\end{equation}
Let us take $p(x)=x^n \in\mathbf{P}_{n}$. From (24), we have
\begin{equation}
\tag{25}
\begin{split}
 C_k&= \frac{(-1)^k2^kk!\delta_{k}}{\pi(2k)!}\int_{-1}^{1}
\biggl(\frac{d^k}{d  x^k }(1-x^2)^{k-1/2}\biggl) x^ndx
\\&=\frac{(-1)^k2^kk!}{\pi(2k)!}\delta_{k}(-1)^k\frac{n!}{(n-k)!}\int_{-1}^{1}
(1-x^2)^{k-1/2} x^{n-k}dx.
\end{split}
\end{equation}
It is easy to show that
\begin{equation}
\tag{26}
\begin{split}
&\int_{-1}^{1}
(1-x^2)^{k-1/2} x^{n-k}dx=\frac{(1+(-1)^{n-k})}{2}\int_{0}^{1}(1-y)^{k-1/2}y^{\frac{n-k+1}{2}-1}dy
\\&=\frac{(1+(-1)^{n-k})}{2}\frac{\Gamma(k+1/2) \Gamma(\frac{n-k+1}{2})}{\Gamma(\frac{k+n+2}{2})}=\frac{(1+(-1)^{n-k})}{2}\frac{(n-k)!(2k)!\pi}{2^{n+k}(\frac{n+k }{2})!(\frac{n-k }{2})!k!}.
\end{split}
\end{equation}
 By (25) and (26), we get
\begin{equation}
\tag{27}
\begin{split}
&C_k= \left\{
\begin{array}{cl}
0, & \text{ if} ~~n-k\equiv1 ~(\text{mod}~2)  \\
\frac{n!\delta_{k}}{2^{n }(\frac{n+k }{2})!(\frac{n-k }{2})! }, & \text{ if} ~~n-k\equiv0 ~(\text{mod}~2). \\
 \end{array}\right.
\end{split}
\end{equation}
From (27), we note that
\begin{equation}
\tag{28}
\begin{split}
x^n=\sum_{k=0}^{n}C_kT_k(x)=\frac{n!}{2^{n-1}}\sum_{\substack{1\leq k\leq n \\ k\equiv1 ~(\text{mod}~2)} } \frac{T_k(x)}{ (\frac{n+k }{2})!(\frac{n-k }{2})! },
\end{split}
\end{equation}
where $n \equiv1~(\text{mod}~2)$.

For $n \equiv0~(\text{mod}~2)$, we have
\begin{equation}
\tag{29}
\begin{split}
x^n= \frac{n!}{2^{n }}\Biggl{\{} \frac{T_0(x)}{\bigl((\frac{n}{2})!\bigl)^2}+2\sum_{\substack{2\leq k\leq n \\ k\equiv0 ~(\text{mod}~2)} } \frac{T_k(x)}{ (\frac{n+k }{2})!(\frac{n-k }{2})! }\Biggl{\}}.
\end{split}
\end{equation}
Let us take $p(x)=B_n(x) \in\mathbf{P}_{n}$. Then
\begin{equation}
\tag{30}
\begin{split}
 C_k&= \frac{(-1)^k2^kk!\delta_{k}}{\pi(2k)!}\int_{-1}^{1}
\biggl( \bigl(\frac{d }{d  x } \bigl)^k(1-x^2)^{k-1/2}\biggl)B_n(x)dx
\\&=\frac{(-1)^k2^kk!\delta_{k}}{\pi(2k)!}(-1)^k\frac{n!}{(n-k)!}\int_{-1}^{1}
(1-x^2)^{k-1/2} B_{n-k}(x)dx
\\&=\frac{ 2^kk!\delta_{k}}{\pi(2k)!} \frac{n!}{(n-k)!}\sum_{l=0}^{n-k}\binom{n-k}{l}B_{n-k-l}\int_{-1}^{1}
(1-x^2)^{k-1/2} x^ldx.
\end{split}
\end{equation}
Now, we compute $\int_{-1}^{1}
(1-x^2)^{k-1/2} x^ldx$.
\begin{equation}
\tag{31}
\begin{split}
\int_{-1}^{1}
(1-x^2)^{k-1/2} x^ldx&=(1+(-1)^l)\int_{0}^{1}
(1-x^2)^{k-1/2} x^ldx
\\&= \left\{
\begin{array}{cl}
0, & \text{ if} ~~l\equiv1 ~(\text{mod}~2)  \\
\frac{l!(2k)!\pi}{2^{2k+l }(\frac{2k+l}{2})!(\frac{l }{2})!k! }, & \text{ if} ~~l\equiv0 ~(\text{mod}~2). \\
 \end{array}\right.
\end{split}
\end{equation}
By (30) and (31), we get
\begin{equation}
\tag{32}
\begin{split}
 C_k&= \frac{ 2^kk!\delta_{k}}{\pi(2k)!}\times \frac{n!}{(n-k)!}\times \frac{(2k)!\pi}{2^{2k}k!}\sum_{\substack{0 \leq l\leq n-k \\ l\equiv0 ~(\text{mod}~2)} }\binom{n-k}{l}B_{n-k-l}\frac{l!}{2^l(\frac{2k+l}{2})!(\frac{l}{2})!}
 \\&=\frac{ n!\delta_{k}}{2^k(n-k)!}\sum_{\substack{0 \leq l\leq n-k \\ l\equiv0 ~(\text{mod}~2)} }\frac{\binom{n-k}{l}B_{n-k-l}l!}{2^l(\frac{2k+l}{2})!(\frac{l}{2})!}.
\end{split}
\end{equation}
Therefore, by (32), we obtain the following theorem.
\begin{theorem}
For $n \in {\mathbb{Z}}_+$, we have
\begin{equation*}
B_n(x)=n!\sum_{0\leq k \leq n} \Biggl( \frac{  \delta_{k}}{2^k(n-k)!}  \sum_{\substack{0 \leq l\leq n-k \\ l\equiv0 ~(\text{mod}~2)} }\frac{\binom{n-k}{l}B_{n-k-l}l!}{2^l(\frac{2k+l}{2})!(\frac{l}{2})!}\Biggl) T_k(x).
\end{equation*}
\end{theorem}
\vspace{3mm}
By the same method, we can derive the following identity:
\begin{equation*}
E_n(x)=n!\sum_{0\leq k \leq n} \Biggl( \frac{  \delta_{k}}{2^k(n-k)!}  \sum_{\substack{0 \leq l\leq n-k \\ l\equiv0 ~(\text{mod}~2)} }\frac{\binom{n-k}{l}E_{n-k-l}l!}{2^l(\frac{2k+l}{2})!(\frac{l}{2})!}\Biggl) T_k(x).
\end{equation*}
Let us take $p(x)=H_n(x) \in\mathbf{P}_{n}$. From (24), we have
\begin{equation}
\tag{33}
\begin{split}
 C_k&= \frac{(-1)^k2^kk!\delta_{k}}{\pi(2k)!}\int_{-1}^{1}
\biggl(  \frac{d^k }{d  x ^k}  (1-x^2)^{k-1/2}\biggl)H_n(x)dx
\\&=\frac{(-1)^k2^kk!\delta_{k}}{(2k)!\pi}\times(-1)^k2^k\frac{n!}{(n-k)!}\int_{-1}^{1}
(1-x^2)^{k-1/2} H_{n-k}(x)dx
\\&=\frac{ 2^{2k}k!\delta_{k}n!}{(2k)!(n-k)!\pi}  \sum_{l=0}^{n-k}\binom{n-k}{l}H_{n-k-l}2^l\int_{-1}^{1}
(1-x^2)^{k-1/2} x^ldx,
\end{split}
\end{equation}
where $H_{n-k-l}$ is the $(n-k-l)$th Hermite number.

By (31) and (33), we get
\begin{equation}
\tag{34}
\begin{split}
 C_k =   n!\delta_{k} \sum_{\substack{0 \leq l\leq n-k \\ l\equiv0 ~(\text{mod}~2)} }\frac{H_{n-k-l} }{(n-k-l)!(\frac{2k+l}{2})!(\frac{l}{2})!}.
\end{split}
\end{equation}
Therefore, by (34), we obtain the following theorem.
\begin{theorem}
For $n \in {\mathbb{Z}}_+$, we have
\begin{equation*}
H_n(x)=n!\sum_{0\leq k \leq n} \Biggl( \delta_{k}  \sum_{\substack{0 \leq l\leq n-k \\ l\equiv0 ~(\text{mod}~2)} }\frac{H_{n-k-l} }{(n-k-l)!(\frac{2k+l}{2})!(\frac{l}{2})!}\Biggl) T_k(x).
\end{equation*}
\end{theorem}
\vspace{3mm}
Let $\mathbf{P}_{n}^{*}=\{p(x)\in  \mathbb{Q}[x] \mid \text{deg}~p(x)\leq n \}$. Then $\mathbf{P}_{n}^{*}$ is an inner product space with
the weighted inner product
  $\langle p(x),~q(x) \rangle =\int_{-1}^{1}\sqrt{1-x^2}p(x)q(x)dx,  \text{ where } ~~   p(x), q (x) \in\mathbf{P}_{n}$ .
Then  $\{U_{0} (x), U_{1} (x), \cdots ,U_{n} (x) \}$ is
 an orthogonal basis for the inner product space  $\mathbf{P}_{n}^{*}$.
 \\
For $p(x) \in\mathbf{P}_{n}^{*}$, let
\begin{equation}
\tag{35}
\begin{split}
p(x)=\sum_{k=0}^{n}C_kU_k(x),
\end{split}
\end{equation}
where
\begin{equation}
\tag{36}
\begin{split}
 C_k &=  \frac{2}{\pi}\langle p(x),~U_k(x) \rangle=\frac{2}{\pi}\int_{-1}^{1}(1-x^2)^{1/2} U_k(x)p(x)dx
 \\&=\frac{(-1)^k2^{k+1}(k+1)!}{(2k+1)!\pi}\int_{-1}^{1}\biggl(\frac{d^k}{dx^k}(1-x^2)^{k+1/2} \biggl)p(x)dx.
\end{split}
\end{equation}
Let us assume that $p(x)=x^n \in\mathbf{P}_{n}^{*}$. Then, by (36), we get
\begin{equation}
\tag{37}
\begin{split}
 C_k &=   \frac{(-1)^k2^{k+1}(k+1)!}{(2k+1)!\pi}\int_{-1}^{1}\biggl(\frac{d^k}{dx^k}(1-x^2)^{k+1/2} \biggl)x^ndx
 \\&=\frac{(-1)^k2^{2k+1}(k+1)!}{(2k+1)!\pi} \times \frac{(-1)^kn!}{(n-k)!}\int_{-1}^{1} (1-x^2)^{k+1/2}  x^{n-k}dx.
\end{split}
\end{equation}
It is easy to show that
\begin{equation}
\tag{38}
\begin{split}
 \int_{-1}^{1}
(1-x^2)^{k+1/2} x^{n-k}dx&=(1+(-1)^{n-k})\int_{0}^{1}(1-x^2)^{k+1/2}x^{n-k}dx
\\&= \left\{
\begin{array}{cl}
0, & \text{ if} ~~n-k\equiv1 ~(\text{mod}~2)  \\
\frac{(n-k)!(2k+2)!\pi}{2^{n+k+2 }(\frac{n+k+2}{2})!(\frac{n-k }{2})!(k+1)! }, & \text{ if} ~~n-k\equiv0 ~(\text{mod}~2). \\
 \end{array}\right.
\end{split}
\end{equation}
Therefore, by (37) and (38), we obtain the following proposition.
\begin{proposition}
For $n \in {\mathbb{Z}}_+$, we have
\begin{equation*}
x^n=\frac{n!}{2^n}   \sum_{\substack{0 \leq k\leq n  \\ k\equiv n ~(\text{mod}~2)} }\frac{k+1}{ (\frac{n+k+2}{2})!(\frac{n-k}{2})!}  U_k(x).
\end{equation*}
\end{proposition}
\vspace{3mm}
Let us consider $p(x)=B_n(x)\in\mathbf{P}_{n}^{*}$. From (36), we have
\begin{equation}
\tag{39}
\begin{split}
 C_k &=   \frac{(-1)^k2^{k+1}(k+1)!}{(2k+1)!\pi}\int_{-1}^{1}\biggl(\frac{d^k}{dx^k}(1-x^2)^{k+1/2} \biggl)B_n(x)dx
 \\&=\frac{(-1)^k2^{ k+1}(k+1)!}{(2k+1)!\pi} \times \frac{(-1)^kn!}{(n-k)!}\int_{-1}^{1} (1-x^2)^{k+1/2}  B_{n-k}(x)dx
  \\&=\frac{ 2^{ k+1}(k+1)!}{(2k+1)!\pi} \times \frac{ n!}{(n-k)!} \sum_{l=0}^{n-k}\binom{n-k}{l}B_{n-k-l}\int_{-1}^{1} (1-x^2)^{k+1/2}  x^ldx.
\end{split}
\end{equation}
It is not difficult to show that
\begin{equation}
\tag{40}
\begin{split}
 \int_{-1}^{1}
(1-x^2)^{k+1/2} x^{l}dx&=(1+(-1)^{l})\int_{0}^{1}(1-x^2)^{k+1/2}x^{l}dx
\\&= \left\{
\begin{array}{cl}
0, & \text{ if } ~~l\equiv1 ~(\text{mod}~2)  \\
\frac{(2k+2)!l!\pi}{2^{2k+2+l }(\frac{2k+2+l}{2})!(k+1)! (\frac{l }{2})!}, & \text{ if } ~~l\equiv0 ~(\text{mod}~2). \\
 \end{array}\right.
\end{split}
\end{equation}
By (39) and (40), we get
\begin{equation}
\tag{41}
\begin{split}
 C_k = \frac{(k+1)n!}{2^k} \sum_{\substack{0 \leq l\leq n-k \\ l\equiv0 ~(\text{mod}~2)} }\frac{B_{n-k-l} }{(n-k-l)!2^l(\frac{2k+l+2}{2})!(\frac{l}{2})!}.
\end{split}
\end{equation}
Therefore, by (41), we obtain the following theorem.
\begin{theorem}
For $n \in {\mathbb{Z}}_+$, we have
\begin{equation*}
B_n(x)=n!\sum_{0\leq k \leq n} \Biggl( \frac{k+1}{2^k} \sum_{\substack{0 \leq l\leq n-k \\ l\equiv0 ~(\text{mod}~2)} }\frac{B_{n-k-l} }{2^l(n-k-l)!(\frac{2k+l+2}{2})!(\frac{l}{2})!}\Biggl) U_k(x).
\end{equation*}
\end{theorem}
\vspace{3mm}
By the same method, we can derive the following identity:
\begin{equation*}
E_n(x)=n!\sum_{0\leq k \leq n} \Biggl( \frac{k+1}{2^k} \sum_{\substack{0 \leq l\leq n-k \\ l\equiv0 ~(\text{mod}~2)} }\frac{E _{n-k-l} }{2^l(n-k-l)!(\frac{2k+l+2}{2})!(\frac{l}{2})!}\Biggl) U_k(x).
\end{equation*}
Let us take $p(x)=H_n(x)\in\mathbf{P}_{n}^{*}$. Then $H_n(x)=\sum_{k=0}^{n}C_kU_k(x)$, with
\begin{equation}
\tag{42}
\begin{split}
 C_k &=   \frac{(-1)^k2^{k+1}(k+1)!}{(2k+1)!\pi}\int_{-1}^{1}\biggl(\frac{d^k}{dx^k}(1-x^2)^{k+1/2} \biggl)H_n(x)dx
 \\&=\frac{ 2^{ 2k+1}(k+1)!n!}{(2k+1)!\pi (n-k)!}   \sum_{l=0}^{n-k}\binom{n-k}{l}2^lH_{n-k-l}\int_{-1}^{1} (1-x^2)^{k+1/2}  x^ldx
  \\&=n!(k+1) \sum_{\substack{0 \leq l\leq n-k \\ l\equiv0 ~(\text{mod}~2)} }\frac{H _{n-k-l} }{ (n-k-l)!}\times \frac{1}{(\frac{2k+l+2}{2})!(\frac{l}{2})!}.
\end{split}
\end{equation}
Thus, by (42) and (43), we get
\begin{equation*}
H_n(x)=n!\sum_{0\leq k \leq n} \Biggl((k+1) \sum_{\substack{0 \leq l\leq n-k \\ l\equiv0 ~(\text{mod}~2)} }\frac{H _{n-k-l} }{ (n-k-l)!(\frac{2k+l+2}{2})!(\frac{l}{2})!}\Biggl) U_k(x).
\end{equation*}

ACKNOWLEDGEMENTS. This research was supported by Basic Science Research Program
through the National Research Foundation of Korea NRF funded by the Ministry of Education,
Science and Technology 2012R1A1A2003786.

%%%%%%%%%%%%%%%%%%%%%%%%%%%%%%%%%%%%%%%%%%%%%%%%%%%%%%%%%%%%%%%%%%%%%%%%%%%%%%%%%%%%%%%%%%%%%%%%%%%%%%%%%%%%%%%%%%%%%%%%%%%%%%%%%%%%%%%%%%%%%%%%%%%%%%%%%%%%%%

\vspace{3mm}

%%%%%%%%%%%%%%%%%%%%%%%%%%%%%%%%%%%%%%%%%%%%%%%%%%%%%%%%%%%%%%%%%%%%%%%%%%%%%%%%%%%%%%%%%%%%%%%%%%%%%%%%%%%%%%%%%%%%%%%%%%%%%%%%%%%%%%%%%%%%%%%%%%%%%%

%
\par\bigskip

\par

\bigskip\bigskip

\begin{center}\begin{large}

{\sc References}

\end{large}\end{center}

\par

\begin{enumerate}

\item[{[1]}] S. Araci, D. Erdal, J. J. Seo, {\it A study on the fermionic $p$-adic $q$-integral representation on $\Bbb Z_p$  associated with weighted $q$-Bernstein and q-Genocchi polynomials},Abstr. Appl. Anal.  Art. ID{\bf649248} (2011), 10 pp.

\item[{[2]}] A. Bayad, T. Kim, {\it  Identities involving values of Bernstein, q-Bernoulli, and q-Euler polynomials}, Russ. J. Math. Phys. {\bf18} (2011), no. 2, 133-143.

\item[{[3]}] L. Carlitz, {\it Note oin the integral of the product of several Bernoulli polynomials},  J.  London Math. Soc.   {\bf34}(1959), 362-363.

\item[{[4]}] M. Can, M. Cenkci, V. Kurt, Y. Simsek, {\it Twisted Dedekind type sums associated with Barnes' type multiple Frobenius-Euler $l$-functions}, Adv. Stud. Contemp. Math. {\bf18} (2009), no. 2, 135-160.

\item[{[5]}]  L. Fox, I. B. Parker, {\it Chebyshev polynomials in numerical analysis},  Oxford University Press, London-New York-Toronto, Ont. 1968.

\item[{[6]}] D. S. Kim, T. Kim, S.-H. Rim, S, H. Lee, {\it Hermite polynomials and their applications associated with Bernoulli and Euler numbers}, Discrete Dynamics in Nature and Society  2012(2012), Article ID {\bf974632}, 13 pp.

\item[{[7]}] D. S. Kim, T. Kim, S.-H. Lee, Y.-H. Kim, {\it Some identities for the product of two Bernoulli and Euler polynomials},Advances in Difference Equations 2012(2012) Article ID {\bf2012:95}, 14 pp.

\item[{[8]}] T. Kim, {\it Identities involving Frobenius-Euler polynomials arising from non-linear differential equations}, J. Number Theory ( Article in Press), 2012.

\item[{[9]}] J. C. Mason,  D.C. Handscomb, {\it  Chebyshev Polynomials}, Chapman  Hall CRC , A CRC Press Company Boca Raton London New York Washington, D.C., 2003.

\item[{[10]}] H. Ozden, I. N. Cangul, Y. Simsek,  {\it Remarks on $q$-Bernoulli numbers associated with Daehee numbers}, Adv. Stud. Contemp. Math. {\bf18} (2009), no. 1, 41-48.

\item[{[11]}] S.-H. Rim, J. Jeong, {\it On the modified $q$-Euler numbers of higher order with weight}, Adv. Stud. Contemp. Math. {\bf22} (2012), no. 1, 93-98.

\item[{[12]}] C. S. Ryoo,  {\it Some identities of the twisted q-Euler numbers and polynomials associated with q-Bernstein polynomials}, Proc. Jangjeon Math. Soc. {\bf14} (2011), no. 2, 239-248.

\item[{[13]}] Y. Simsek, {\it Generating functions of the twisted Bernoulli numbers and polynomials associated with their interpolation functions}, Adv. Stud. Contemp. Math. {\bf16} (2008), no. 2, 251-278.

\item[{[14]}] Y. Simsek, {\it  Theorems on twisted $L$-function and twisted Bernoulli numbers}, Adv. Stud. Contemp. Math. {\bf11} (2005), no. 2, 205-218.

\item[{[15]}] Y. Simsek, M. Acikgoz, {\it A new generating function of $(q-)$Bernstein-type polynomials and their interpolation function}, Abstr. Appl. Anal. 2010, Art. ID {\bf769095}, 12 pp.

\end{enumerate}

\par

\bigskip\bigskip

\par

\bigskip\bigskip

\par\noindent

\mpn { \bpn {\small Dae San {\sc Kim} \mpn Department of
Mathematics, \pn Sogang University,  Seoul 121-742, Republic of Korea 

\bpn {\small TaeKyun {\sc Kim} \mpn Department of Mathematics,
\pn Kwangwoon University, Seoul 139-701, Republic of Korea \pn {\it E-mail:}\ {\sf tkkim@kw.ac.kr, taekyun64@hotmail.com} }

\bpn {\small Sang-Hun  {\sc Lee} \mpn  Division of General Education,\pn Kwangwoon University, Seoul 139-701, Republic of Korea 

\end{document}